\documentclass[16pt]{amsart}
\usepackage{amssymb,amsmath,amsthm,setspace,verbatim}
\onehalfspace

\newtheorem{thm}{Theorem}
\newtheorem*{thm*}{Sturm's Theorem}

\theoremstyle{definition}

\theoremstyle{remark}

\begin{document}

\title{Root Configurations of Real Univariate Cubics and Quartics}

\author{Elias Gonzalez}
\address{San Antonio, Texas}
\email{elias.gonzalez@nisd.net}

\author{David A. Weinberg}
\address{Department of Mathematics and Statistics, Texas Tech University, 
Lubbock, TX 79409}
\email{david.weinberg@ttu.edu}
\urladdr{www.math.ttu.edu/~dweinber}

\begin{abstract}
For the general monic cubic and quartic with real coefficients, polynomial conditions on the coefficients are derived as directly and as simply as possible from the Sturm sequence that will determine the real and complex root multiplicities together with the order of the real roots with respect to multiplicity.
\end{abstract}

\maketitle

\section{Introduction}

Everyone knows about the discriminant $b^{2}-4ac$ of the quadratic $ax^{2}+bx+c$.  Any univariate polynomial $f(x)$ has a discriminant, which is (essentially) the resultant of $f(x)$ and its derivative.  The discriminant is the determinant of a certain matrix formed from the coefficients of $f(x)$ and its derivative.  The significance of the discriminant is essentially that it vanishes if and only if the polynomial $f(x)$ has a multiple root.  By looking at determinants of submatrices of matrices giving resultants, one can get more information about the multiple roots.  In fact, resultants were used to get the complex root multiplicities in [6].  In the 1990's, other authors studied root multiplicities and obtained extensive results.  In [4], [7], and [8], those authors developed and applied the notions of complete discrimination system, multiple factor sequence, and revised sign list, and proved general theorems about conditions for root multiplicities.\\

In this paper, more detailed information about the roots of univariate polynomials (e.g. cubics and quartics) with real coefficients will be obtained by using a variant of the Euclidean algorithm, called the Sturm sequence.  We want to find polynomial conditions in terms of the coefficients that will determine not only the real and complex root multiplicities, but also the relative position on the real line of the real roots with respect to their multiplicities.  These results are bound to have applications in the future since one of the major paradigms of real algebraic geometry is to regard a multivariable polynomial as a polynomial in one variable with parameters.  For example, they can be applied to the study of an affine hypersurface regarded as a branched cover of the hyperplane obtained by setting the last coordinate equal to zero.  The significance of this paper is twofold: (1) The polynomial conditions on the coefficients that determine the order of the real roots with regard to multiplicity are new, and (2) The method of proof, based on the single unifying idea of the Sturm sequence, is the simplest and most transparent possible when considering the specific cases of cubics and quartics.\\

\quad Consider the general cubic $ax^{3}+bx^{2}+cx+d$.  Then without loss of generality, we can instead consider the monic cubic obtained by dividing the coefficients by $a$.  This monic cubic has the form $x^{3}+px^{2}+qx+r$, where $p=\frac{b}{a}$, $q=\frac{c}{a}$, and $r=\frac{d}{a}$.  Consideration of the monic cubic will make various computations simpler and free us from considering the cases where $a=0$.  In general, experience has taught us that it is best to consider the space of general monic polynomials (including the term of degree $n-1$).  Here are the possible root configurations for the cubic:\\
\flushleft
\begin{enumerate}
\item 3 single real roots
\item 1 single real root and 2 complex conjugate roots
\item 1 double real root and 1 single real root
\begin{enumerate}
\item double real root $<$ single real root
\item single real root $<$ double real root
\end{enumerate}
\item 1 triple real root\newline
\end{enumerate}

\quad We want to find polynomial conditions on $p$, $q$, and $r$ that will determine which of these possible configurations hold.  We will use the Sturm sequence.  The novelty here is that the calculations will be done with symbolic coefficients, in other words, over the function field of the coefficients.  At each stage, the polynomials will have coefficients that are rational functions of the coefficients of the original polynomial $f(x)$.  For degree four and higher, such calculations are too tedious and lengthy to do by hand, but can be done rapidly and conveniently by using Maple software.\\
\quad As all undergraduate mathematics majors know, the Euclidean algorithm gives the greatest common divisor of two univariate polynomials $f(x)$ and $g(x)$ with coefficients in a field.  When applied to $f(x)$ and $f'(x)$, it gives the product of the multiple roots of $f(x)$, each counted with multiplicity one less than its multiplicity as a root of $f(x)$.\\
\quad If $c=\{c_{1},c_{2},$\ldots,$c_{m}\}$ is a finite sequence of real numbers, then the \emph{number of variations in sign} of $c$ is defined to be the number of $i$, $1\leq i\leq m-1$ such that $c_{i}c_{i+1}<0$, after dropping the 0's in $c$.  The \emph{Sturm sequence} for $f(x)$ is defined to be $f_{0}(x)=f(x)$, $f_{1}(x)=f'(x)$, $f_{2}(x)$, \ldots, $f_{s}(x)$, where\\
\newpage
\quad $f_{0}(x)=q_{1}(x)f_{1}(x)-f_{2}(x)$ \quad deg$f_{2}(x)<$deg$f_{1}(x)$\\
\quad \vdots \\
\quad $f_{i-1}(x)=q_{i}(x)f_{i}(x)-f_{i+1}(x)$ \quad deg$f_{i+1}(x)<$deg$f_{i}(x)$\\
\quad \vdots \\
\quad $f_{s-1}(x)=q_{s}(x)f_{s}(x)$ \quad $(f_{s+1}(x)=0)$\\
In other words, perform the Euclidean algorithm and change the sign of the remainder at each stage.  It should be noted here that there exists a more general definition of Sturm sequence, but the one given here suits our purpose.\\

\begin{thm*} Let $f(x)$ be a polynomial of positive degree in $\mathbb{R}[x]$ and let $f_{0}(x)=f(x)$, $f_{1}(x)=f'(x)$, \ldots, $f_{s}(x)$ be the Sturm sequence for $f(x)$.  Assume $[a,b]$ is an interval such that $f(a)\neq0$ and $f(b)\neq0$.  Then the number of distinct real roots of $f(x)$ in $(a,b)$ is $V_{a}-V_{b}$, where $V_{c}$ denotes the number of variations of sign of $\{f_{0}(c)$, $f_{1}(c)$, \ldots, $f_{s}(c)\}$.\end{thm*}

Proofs can be found in [1], [2], [3], and [5].  We can get the total number of real roots by looking at the limits as $a\to-\infty$ and $b\to+\infty$.  Thus, the total number of distinct real roots will depend only on the leading terms of the polynomials in the Sturm sequence. In the rest of the paper, when we refer to a use of Sturm's Theorem, or if we say 'Sturm analysis', it is left to the reader to count the sign changes at $-\infty$ and $\infty$ for all of the sign possibilities of the leading coefficients.\\
We will also use the following theorem, which the reader can prove as an exercise(or see [1], page 102):

\begin{thm} If the discriminant of a real polynomial is not zero, then the sign of the discriminant is $(-1)^{r}$, where $r$ is the number of complex conjugate pairs of roots.\end{thm}

\quad Now consider the cubic, $f(x)=x^{3}+px^{2}+qx+r$.  Calculate the Sturm sequence for $f(x)$.  We obtain polynomials of degree one and degree zero in $x$ with coefficients that are rational functions of $p$, $q$, and $r$, which we will name $gcddeg1$ and $gcddeg0$, respectively.\\
\quad If $gcddeg0\neq0$, then there are three distinct roots.  (It turns out, not surprisingly, that this condition is equivalent to the discriminant not being zero.)  The cubic is so simple that the number of real roots and the number of complex roots is determined by Theorem 1. (For the quartic, we will need to use both Theorem 1 and Sturm's Theorem.)\\
\quad If $gcddeg0=0$ and the leading coefficient of $gcddeg1\neq0$, this means the cubic has a double root, which is, in fact, the zero of $gcddeg1$, and so one can solve for it.  In this case, we know there is one double real root and one single real root.  Which one is bigger?\\
\quad To find conditions which will determine which root is larger, solve $gcddeg1=0$ for the double root.  Then, make a change of variable, $x=y+doubleroot$, so as to translate the double root to the origin.  After this translation, $y^{2}$ becomes a factor of the resulting cubic in $y$, and one can now solve for the single root.  Since translation does not affect the relative position of the roots, if the new single root is positive, then it is bigger than the double root, and if the new single root is negative, then it is less than the double root, and the same holds for the polynomial in $x$.\\
\quad If $gcddeg0=0$ and $gcddeg1=0$, then there is a triple root.\\

\quad The unifying idea of the Sturm sequence is the optimal and most transparent approach to the classical questions about roots of univariate polynomials answered in this paper.  The computer is used to perform long divisions, calculate discriminants, and perform factoring, all of which are too tedious to do by hand, at least in the case of the quartic.  The novelty here is to study the Sturm sequence in complete detail for a general monic polynomial with variable coefficients.  Let us note here that one has to keep careful track of whether or not certain coefficients of the polynomials in the Sturm sequence are zero.  This leads to many special cases to consider, the number of which increases rapidly with the degree of $f(x)$.  In fact, if $f(x)$ has degree $n\geq3$, then the number of special cases is $2^{n-2}-1$.\\

\quad In the following sections, we will use some suggestive notation. $gcddeg n$ denotes the polynomial of degree $n$ in the Sturm sequence, and $new m gcddeg n$ denotes the polynomial of degree $n$ in the Sturm sequence for the $m^{\mathrm{th}}$ special case.\\

\quad The derivation and proofs in this paper establish a partition of the spaces of monic cubic and quartic polynomials with real coefficients according to the real and complex root multiplicities, taking into account the order of the real roots with regard to multiplicity.  Two other, coarser, results follow immediately: a partition of the same spaces by real and complex root multiplicities (without order) and a partition of the space of monic complex cubics and quartics by complex root multiplicities.\\
\newpage

\section{Cubics}
\quad The Sturm sequence for the cubic is: \\
\begin{center}
$x^{3}+px^{2}+qx+r$\\
$3x^{2}+2px+q$\\
$gcddeg1 = \frac{2}{9}(p^{2}-3q)x-r+\frac{pq}{9}$\\
$gcddeg0 = \frac{9(-4q^{3}+p^{2}q^{2}+18rpq-4rp^{3}-27r^{2})}{4(p^{2}-3q)^{2}}$\\
\end{center}

Let us remark that the discriminant of the cubic is $D = -4q^{3}+p^{2}q^{2}+18rpq-4rp^{3}-27r^{2}$.\\

\subsection{Root multiplicities.}

\subsubsection{3 distinct roots.} By Theorem 1, if $D > 0$, then there are 3 real single roots, and if $D < 0$, then there is 1 real root and 2 complex conjugate roots.\\

\subsubsection{1 double root and 1 single root.}

\quad In order to have 1 double root and 1 single root, we must have $D = 0$ and  $p^{2}-3q \neq 0$. \\



\subsubsection{triple root.} The conditions for a triple root are $D = 0$ and $p^2-3q = 0$.  (Note that if $3p^2-q = 0$, then the constant term of $gcddeg1$, which is $(pq/9)-r$, becomes $(p^3/27)-r)$, and the discriminant of the cubic becomes $-(p^3-27r)^2/27$, and since the discriminant is zero, the constant term of $gcddeg1$ must also be zero; therefore we must have a triple root.)\\


\subsection{Order of the double and single roots.}  If $D=0$ and $p^{2}-3q\neq0$, then the greatest common divisor of the cubic and its derivative has degree one, and it follows from Sturm's Theorem that $p^2-3q > 0$.  But then the root of that polynomial is the double root.  Solving for $x$ in gcddeg1, we get\\

\begin{equation*}
		\mbox{ doubleroot } = -\left(\frac{-9r+pq}{2(p^{2}-3q)}\right).
\end{equation*}

We will now translate this double root to the origin by making a simple substitution $x=y+$doubleroot.  This will not change the relative position of the single and double roots. The resulting polynomial in $y$ will then have a factor of $y^{2}$ (since the double root is at the origin), and we will be able to solve for the single root.  (It is interesting to note that if $D=0$, we can solve the cubic.)\\

By substituting $x=y+$doubleroot, we get \\
\begin{eqnarray*} y^{3}+\left(-\frac{3(9r-pq)}{2(3q-p^{2})}+p\right)y^{2}&+&\left[\frac{3(9r-pq)^{2}}{4(3q-p^{2})^{2}}+ q-\frac{p(9r-pq)}{3q-p^{2}}\right]y-\frac{(9r-pq)^{3}}{8(3q-p^{2})^{3}}\\  
&+& r - \frac{q(9r-pq)}{2(3q-p^{2})}+\frac{p(9r-pq)^{2}}{4(3q-p^{2})^{2}}.
\end{eqnarray*}

We know that the coefficient of $y$ and the constant term have to be zero, but let us show that here.  Simplifying, the constant term is 
\begin{equation*}
  \frac{-(729r^{3}-729r^{2}pq+135rp^{2}q^{2} + 	   17p^{3}q^{3}+108rq^{3}-72rp^{4}q+8rp^{6}-36pq^{4})}{8(3q-p^{2})^{3}}
\end{equation*}

\begin{equation*} \frac{-(-2p^{5}q^{2}+162r^{2}p^{3})}{8(3q-p^{2})^{3}}
= \frac{-(27r+2p^{3}-9pq)(4q^{3}-p^{2}q^{2}-18rpq+4rp^{3}+27r^{2})}{8(3q-p^{2})^{3}}.
\end{equation*} 
Since $D$ is a factor of the numerator, this is zero.  Simplifying, the coefficient of $y$ is $\frac{9(4q^{3}-p^{2}q^{2}-18rpq+4rp^{3}+27r^{2})}{4(3q-p^{2})^{2}}$.  Once again, since $D$ is a factor of the numerator, this is zero.  So, the single root for the polynomial in $y$ is $\frac{27r-9pq+2p^{3}}{2(p^{2}-3q)}$.  Since the substitution did not change the relative position of the roots, and since we know that $p^2-3q > 0$, we know that if $27r-9pq+2p^{3} > 0$, then the single root is greater then the double root, and if $27r-9pq+2p^{3} < 0$, then the single root is less than the double root.\\

To summarize,\\

\textbf{Real Root Configurations.}
\begin{center}
\begin{tabular}{l r}
1. 3 distinct real roots & $D > 0$\\
2. 1 real root and 2 complex conjugate roots & $D < 0$\\
3. 1 double root and 1 single root & $D=0$ and $p^{2}-3q\neq0$\\
3a. single root $>$ double root & $27r-9pq+2p^{3} > 0$\\
3b. single root $<$ double root & $27r-9pq+2p^{3} < 0$\\
4. triple root & $D=0$ and $p^{2}-3q=0$\newline
\end{tabular}
\end{center}

\par \quad An immediate corollary of the above is the following:\newline

\textbf{Complex Root Configurations.}
\begin{center}
\begin{tabular}{l r}
1. 3 distinct roots & $D\neq0$\\
2. 1 double root and 1 single root & $D=0$ and $p^{2}-3q\neq0$\\
3. triple root &  $D=0$ and $p^{2}-3q=0$\\
\end{tabular}
\end{center}
\newpage

\section{Quartics}
\quad Consider the partition of the space of monic quartics, $x^{4}+px^{3}+qx^{2}+rx+s$, according to the following root configurations:
\begin{enumerate}
\item 4 distinct real roots
\item 2 distinct real roots and 2 distinct complex conjugate roots
\item 4 distinct complex roots
\item 1 double real root and two distinct real roots
\begin{enumerate}
\item single root $<$ double root $<$ single root
\item double root $<$ single root $<$ single root
\item single root $<$ single root $<$ double root
\end{enumerate}
\item 1 double real root and 2 distinct complex conjugate roots
\item 2 real double roots
\item 2 complex conjugate double roots
\item 1 triple root and 1 single root
\begin{enumerate}
\item triple root $<$ single root
\item single root $<$ triple root
\end{enumerate}
\item 1 quadruple root
\end{enumerate}
We will find polynomial conditions on the coefficients $p$, $q$, $r$, and $s$ that will determine to which of these classes the polynomial belongs. 

\quad We use the computer software Maple to calculate the Sturm sequence for the quartic.  One obtains remainder polynomials of degrees 2, 1, and 0, which we suggestively name $gcddeg2$, $gcddeg1$, and $gcddeg0$, respectively.  Here is a summary of the generic analysis of the general quartic that will be carried out in this section:\\

\quad If $gcddeg0 \neq0$, then there are 4 distinct single roots.  Then perform a Sturm sequence on the quartic from $-\infty$ to $\infty$, as defined in the introduction.  The signs of the leading coefficients will yield the number of real and complex roots.\\

\quad If $gcddeg0=0$ and the leading coefficient of $gcddeg1\neq0$, then this means the quartic has a real double root and two single roots.  The Sturm sequence determines whether the two single roots are both real or both complex.  This is seen only to depend on the numerator of the leading coefficient of $gcddeg1$.  If the two single roots are real, then we want to know the relative position of the double root and the two single roots.  The double root is the zero of $gcddeg1$.  Now solve for this double root.  Translate the double root to the origin via $x = y + doubleroot$. Now, $y^{2}$ is a factor of the resulting quartic, and there is a leftover quadratic polynomial in $y$.  Now if the constant term of this leftover quadratic, which is the product of the roots, is negative, then the double root lies in between the two single roots.  Next, if the constant term is positive, then the two single roots lie on the same side of the double root, determined by the coefficient of $y$, which is minus the sum of the roots.  This will exactly determine the relative position of the real roots, since the double root has been placed at the origin.\\

\quad If $gcddeg0 = 0$, and $gcddeg1 = 0$, and the leading coefficient of $gcddeg2 \neq 0$, then there are either two double roots or one triple root and one single root.  The Sturm sequence determines whether there are two or zero real roots.  If there are zero real roots, then there are two complex conjugate double roots.  If there are two real roots, then more work is required to distinguish between the cases of two real double roots and one triple real root and one single real root.\\

\quad If $gcddeg0$, $gcddeg1$, and $gcddeg2$ are all equal to zero, then there is a quadruple root.\\

\quad The Sturm sequence for the quartic is:
\begin{center}
$x^{4}+px^{3}+qx^{2}+rx+s$\\
$4x^{3}+3px^{2}+2qx+r$\\
$gcddeg2=\frac{1}{16}(3p^{2}-8q)x^{2}+(\frac{1}{8}pq-\frac{3}{4})x+\frac{1}{16}pr-s$\\
$gcddeg1=\frac{32(-4q^{3}+p^{2}q^{2}+16qs+14rpq-6p^{2}s-3rp^{3}-18r^{2})x}{(8q-3p^{2})^{2}}-\frac{16(-3pr^{2}+4q^{2}+9sp^{3}-p^{2}rq-32spq+48sr)}{(8q-3p^{2})^{2}}$\\
$gcddeg0=$\\
$\frac{1}{(64(-p^{2}q^{2}+3rp^{3}+6p^{2}+4q^{3}-14pqr-16sq+18r^{2})^{2})}(1584p^{4}r^{2}sq-7296p^{2}r^{2}q^{2}s-5120prq^{4}s-12288prs^{2}q^{2}-1584p^{5}rq^{2}s+4992p^{3}rq^{3}s+9216p^{3}rs^{2}q+162p^{7}qrs+16384s^{3}q^{2}+1024sq^{6}-256r^{2}q^{5}-1728r^{4}q^{2}-243p^{8}s^{2}-36p^{7}r^{3}-243p^{4}r^{4}+2304p^{4}s^{3}-8192q^{4}s^{2}-1728p^{5}rs^{2}+1152pr^{3}q^{3}+9p^{6}r^{2}q^{2}-84p^{4}r^{2}q^{3}-1120p^{3}r^{3}q^{2}-54p^{6}r^{2}s+354p^{5}r^{3}q+256p^{2}r^{2}q^{4}+1296p^{2}r^{4}q+336sp^{4}q^{4}-1024sp^{2}q^{5}+9216sq^{3}r^{2}-9792p^{4}q^{2}s^{2}+15360p^{2}q^{3}s^{2}-12288p^{2}s^{3}q-36p^{6}q^{3}s+2592p^{6}qs^{2})$\\
\end{center}

The discriminant of the quartic is $D=-80q^{2}srp+18qsrp^{3}-4q^{3}p^{2}s+p^{2}q^{2}r^{2}+144s^{2}qp^{2}+144sqr^{2}+18r^{3}pq-6p^{2}sr^{2}+16q^{4}s-4q^{3}r^{2}-128s^{2}q^{2}-27p^{4}s^{2}-4r^{3}p^{3}-27r^{4}+256s^{3}-192prs^{2}$.
Notice that the numerator of $gcddeg0$ is $(8q-3p^{2})^{2}D$.
Let $L1$ denote the numerator of the coefficient of $x$ in $gcddeg1$ (without the 32).\\

\subsection{Root multiplicities.}
\subsubsection{4 distinct roots.} Let $L2$ denote the appropriate factor of the leading coefficient of $gcddeg2$.  By Theorem 1, if $D < 0$, then the quartic has 2 real roots and 2 complex conjugate roots.  By Sturm's Theorem, if $D > 0$ and $L2 > 0$ and $L1 >0$, then the quartic has 4 distinct real roots.  By Sturm's Theorem and Theorem 1, if $D > 0$ and ($L2 \leq 0$ or $L1 \leq 0$), then the quartic has no real roots and 2 distinct pairs of complex conjugate roots.\\

\subsubsection{1 double root and 2 single roots.} We have that $D = 0$ and the gcd has degree 1.  If $L2 \neq 0$, then by Sturm's Theorem, if $L1 > 0$, there are 3 real roots, hence one double real root and 2 single real roots, while if $L1 < 0$, then there is 1 real root, hence the 2 single roots are complex conjugate.  In the special case that $L2 = 0$, solving for $q$ gives $q = 3p^2/8$.  In this case it is not possible to get 3 real roots.  Plugging $q = 3p^2/8$ into $L1$ gives $-9(-16r+p^3)^2/128$, and so we see that in this case $L1$ is forced to be negative.

\subsubsection{2 double roots or 1 triple root and 1 single root.}  We want the gcd to have degree 2, so we have $D = 0$ and $L1 = 0$ and $L2 \neq 0$ (by definition of the Sturm sequence, the constant term of $gcddeg1$ is forced to be zero). The Sturm analysis shows that if $3p^2-8q > 0$, then the two roots are real, while if $3p^2-8q < 0$, the two roots are complex.  The discriminant of $gcddeg2$, which will be denoted $D6$, will give all the information about these roots.\\
 $D6 = -{\frac {3}{64}}\,r{p}^{3}+3/4\,{p}^{2}s+{\frac {1}{64}}\,{p}^{2}{q}^{
2}-1/16\,rqp+{\frac {9}{16}}\,{r}^{2}-2\,sq$\\
If $D6 > 0$, there are 2 double real roots, if $D6 = 0$, there is one triple real root and 1 single real root, and if $D6 < 0$, there is one pair of complex conjugate double roots.  However, this will not be our final answer.  (It should be noted here that if one solves $L1 = 0$ for $s$ and then recomputes the Sturm sequence, one can obtain "simpler" conditions for the distinction between triple - single and double - double real, if one combines that computation with a formula for the triple root, which will be derived below.)\\

\subsubsection{quadruple root.} The conditions for a quadruple root are $D = 0$ and $L1 = 0$ and $L2 = 0$.  Computer algebra shows that this is equivalent to setting all the coefficients of $gcddeg2$ equal to zero.\\

\subsection{Order of the real roots with respect to multiplicity.}

\subsubsection{1 double and 2 single roots.} Solve $gcddeg1 = 0$ for the double root.  The result is\\
\begin{center}
$doubleroot=\frac{-3pr^{2}+4q^{2}r+9sp^{3}-p^{2}rq-32spq+48sr}{2(-4q^{3}+p^{2}q^{2}+16qs+14rpq-6p^{2}s-3rp^{3}-18r^{2})}=C_0/2L1$, where\\
\end{center}
\begin{center}
$C_0=48sr-32spq-p^2rq+9sp^3+4q^2-3pr^2$ \\
\end{center}
We will now translate this double root to the origin by substituting $x=y+doubleroot$ into the original quartic.  $y^{2}$ will be a factor of the resulting quartic. Computer algebra shows that the remaining quadratic factor is\\
\begin{center}
$leftoverquadratic=y^{2}+[p+\frac{2(-3pr^{2}+4q^{2}r+9sp^{3}-p^{2}rq-32spq+48sr)}{-4q^{3}+p^{2}q^{2}+16qs+14rpq-6p^{2}s-3rp^{3}-18r^{2}}]y+\frac{3p(-3pr^{2}+4q^{2}r+9sp^{3}-p^{2}rq-32spq+48sr)}{2(-4q^{3}+p^{2}q^{2}+16qs+14rpq-6p^{2}s-3rp^{3}-18r^{2})}+q+\frac{3(-3pr^{2}+4q^{2}r+9sp^{3}-p^{2}rq-32spq+48sr)^{2}}{2(-4q^{3}+p^{2}q^{2}+16qs+14rpq-6p^{2}s-3rp^{3}-18r^{2})^{2}}=y^2+(p+2C_0/L1)y+(3/2)pC_0/L1+q+(3/2)C_0^2/(L1)^2$.\\
\end{center}


If the constant term of $leftoverquadratic$ is negative, then we have $single < double < single$.  If the constant term is positive and the coefficient of $y$ in $leftoverquadratic$ is negative, then we have $double < single < single$, while if the coefficient of $y$ is positive, then we have $single < single < double$.  (It is not necessary to consider the special case $L2 = 0$ because in that case the Sturm analysis shows that the 2 single roots cannot be real.)\\

Observe that the constant term of $leftoverquadratic$ is
\begin{center}
$(3/2)p(C_0/L1)+q+(3/2)C_0^2/L1^2 = $
\end{center}
\begin{center}
$\frac{3C_0^2+2qL1^2+3pC_0L1}{2L1^2}$
\end{center}
(and the sign depends only on the numerator).

Observe that the coefficient of $y$ in $leftoverquadratic$ is
\begin{center}
$p+2C_0/L1 = \frac{pL1+2C_0}{L1} = \frac{(p^3-4pq+8r)(3rp-12s-q^2)}{-L1}$.
\end{center}

(Note that Maple factored the numerator.)\\
 
  Let $D2 = 3C_0^2 + 2q(L1)^2 + 3pC_0L1$, the numerator of the constant term in $leftoverquadratic$. Since we are in the case of one real double root and two real single roots, we know that $L1>0$.  Define $D3 = (p^3-4pq+8r)(q^2-3pr+12s)$.  If $D2 < 0$, then we have $single < double < single$.  If $D2 > 0$ and $D3 > 0$, then we have $single < single < double$.  If $D2 > 0$ and $D3 < 0$, then we have $double < single < single$.

\subsubsection{1 triple and 1 single root.}  We have $D = 0$ and $L1 = 0$.  Since we have a triple root, $gcddeg2$ will have a double root, which is $-1/2$ times the coefficient of $x$, after dividing by the leading coefficient of $gcddeg2$.  Thus
\begin{center}
$tripleroot = \frac{6r-pq}{3p^2-8q}$
\end{center}

Now translate the triple root to the origin by substituting $x = y + tripleroot$ into the original quartic.  The result is 
\begin{center}
$y^4+(p+\frac{4(6r-pq}{3p^2-8q})y^3$
\end{center}

We can now solve for the single root:
\begin{center}
$singleroot = \frac{3(4pq-8r-p^3)}{3p^2-8q}$
\end{center}

Since we know that in this case $3p^2-8q > 0$, the relative location of the single root is determined by the sign of the numerator.  Thus, if $4pq-8r-p^3 > 0$, we have $triple < single$, while if $4pq-8r-p^3 < 0$, we have $single < triple$.

\subsubsection{1 triple and 1 single vs. 2 double real roots reconsidered} If we solve $L1 = 0$ for $s$ and recompute the Sturm sequence, then we can obtain some interesting information.  (This cannot be conveniently done for the quintic since the numerator of the leading coefficient of $gcddeg1$ in the case of a quintic does not having any variable appearing linearly.) Solve $L1 = -4q^{3}+p^{2}q^{2}+16qs+14rpq-6p^{2}s-3rp^{3}-18r^{2} = 0$ for $s$ and plug into the original quartic.  This gives\\

\begin{center}
$newquartic1=x^{4}+px^{3}+qx^{2}+rx-\frac{4q^{3}-p^{2}q^{2}+18r^{2}-14rpq+3rp^{3}}{2(3p^{2}-8q)}$.\\
\end{center}

The Sturm sequence of $newquartic1$ yields the following remainders:\\
\begin{center}
$new1gcddeg2 = (\frac{3}{16}p^{2}-\frac{1}{2}q)x^{2}+(\frac{1}{8}pq-\frac{3}{4}r)x-\frac{32q^{3}-8p^{2}q^{2}+144r^{2}-120rpq+27rp^{3}}{16(3p^{2}-8q)}$\\
$new1gcddeg0 = \frac{8(256rq^{3}+360rp^{2}q^{2}-216rp^{4}q+27rp^{6}-128pq^{4}+68p^{3}q^{3}-1296r^{2}pq-9p^{5}q^{2}+324r^{2}p^{3}+864r^{3})}{(3p^{2}-8q)^{3}}$\\
\end{center}

Let the numerator of $new1gcddeg0$ (without the factor of 8) be $D4$.  Maple computer algebra yields $D4=\left( 8\,r+{p}^{3}-4\,pq \right)  \left( 27\,r{p}^{3}-9\,{p}^{2}{q}^
{2}-108\,pqr+32\,{q}^{3}+108\,{r}^{2} \right)$.  Observe that the discriminant $D$ of $newquartic1$ is $-{\frac { \left( 8\,r+{p}^{3}-4\,pq \right) ^{2} \left( 27\,r{p}^
{3}-9\,{p}^{2}{q}^{2}-108\,pqr+32\,{q}^{3}+108\,{r}^{2} \right) ^{2}}{4 \left( 3\,{p}^{2}-8\,q \right) ^{3}}}$.\\

\quad If $D4 = 0$, then $newquartic1$ has either two double roots or a triple root and a single root.  The Sturm analysis of $newquartic1$ determines whether there are two real roots or zero real roots.  If $3p^{2}-8q<0$, then there are zero real roots, i.e. two complex conjugate double roots.  If $3p^{2}-8q>0$, then there are two real roots, i.e. two real double roots or one real triple root and one real single root.  It is now necessary to find the conditions that will distinguish between two real double roots and one triple and one single root.\\

\quad Consider the the discriminant of $new1gcddeg2$, which will be denoted $D5 = -\frac{27}{16}r^{2}-\frac{1}{2}q^{3}+\frac{27}{16}rpq+\frac{9}{64}p^{2}q^{2}-\frac{27}{64}rp^{3}$.  If $D5=0$, $new1gcddeg2$ has a double root, which is $-\frac{1}{2}$ times the coefficient of $x$.  In this case, $newquartic1$ has a real triple root and a real single root; translate the triple root to the origin.  Thus,\\
\begin{center}
$new1tripleroot=\frac{-pq+6r}{3p^{2}-8q}$
\end{center}

Plug $x = y + new1tripleroot$ into $newquartic1$.  Since we have shifted the triple root to zero, we know the resulting quartic will have coefficients of zero for $y^{2}$ and $y$, and the constant term is zero.  The resulting quartic is $y^{4}+(p+\frac{4(-pq+6r)}{3p^{2}-8q})y^{3}$. If we factor out $y^{3}$, we are left with a linear equation, which we may solve for the single root of the shifted $newquartic1$.
\begin{center}
$new1singleroot =-\,{\frac {3(8\,r+{p}^{3}-4\,pq)}{3\,{p}^{2}-8\,q}}$.\\
\end{center}

Recall that in order to have two real roots, there must be the condition $3p^{2}-8q>0$.  If $p^{3}-4pq+8r<0$, then the single root is greater than the triple root, and if $p^{3}-4pq+8r>0$, then the single root is less than the triple root.  Note that $p^{3}-4pq+8r\neq0$ because the triple root and single root are both distinct, and the triple root of the translated quartic is at zero.\\

\quad If $D5\neq0$, then there are two double real roots.  Notice that $D4=-64(p^{3}-4pq+8r)(D5)$.  Since we are in the case where $D4=0$, this forces $p^{3}-4pq+8r=0$.  Thus, we have shown that the conditions for triple and single root are $D=0$ and $L1=0$ and $3p^{2}-8q>0$ and $p^{3}-4pq+8r\neq0$, and the conditions for two real double roots are $D=0$ and $L1=0$ and $3p^{2}-8q>0$ and $p^{3}-4pq+8r=0$.\bigskip
\quad To summarize,\\

\textbf{Notation.}
\begin{multline*}
D=18\,{p}^{3}rqs-4\,{q}^{3}{p}^{2}s+144\,s{r}^{2}q+{q}^{2}{p}^{2}{r}^{2}-192\,pr{s}^{2}+144\,q{p}^{2}{s}^{2}+18\,p{r}^{3}q-4\,{p}^{3}{r}^{3}-128\,{q}^{2}{s}^{2}\\
+16\,{q}^{4}s-4\,{q}^{3}{r}^{2}-27\,{p}^{4}{s}^{2}-80\,pr{q}^{2}s-6\,{p}^{2}{r}^{2}s+256\,{s}^{3}-27\,{r}^{4}; \text{ the discriminant of the quartic}\\
L1=p^{2}q^{2}-3rp^{3}-6p^{2}s-4q^{3}+14pqr+16sq-18r^{2}; \frac{1}{32} \text{ times the numerator of } gcddeg1\\
C_0=48sr-32spq-p^2rq+9sp^3+4q^2-3pr^2\\
D2=3C_0^2 + 2q(L1)^2 + 3pC_0L1;\, \text{the numerator of the constant term in}\   leftoverquadratic.\\
D3=(8r-4pq+p^{3})(q^{2}-3pr+12s)\\
\end{multline*}

\textbf{Real Root Configurations.}
\begin{tabular}{l r}
1. 4 distinct real roots & $3p^{2}-8q > 0$ and $L1 > 0$ and $D>0$\\
2. 2 real roots and 2 complex conjugate roots & $D<0$\\
3. 4 complex single roots & $D>0$ and ($3p^{2}-8q\leq0$ or $L1\leq0$)\\
4. 1 real double root and & $D=0$ and $L1>0$\\ \quad2 real single roots & \\
4a. single root $<$ double root $<$ single root & $D2<0$\\
4b. double root $<$ single root $<$ single root & $D2>0$ and $D3<0$\\
4c. single root $<$ single root $<$ double root & $D2>0$ and $D3>0$\\
5. 1 real double root and & $D=0$ and $L1<0$\\ \quad2 complex conjugate roots & \\
6. 2 real double roots & $D=0$ and $L1=0$ and\\ & $3p^{2}-8q>0$ and $p^{3}-4pq+8r=0$\\
7. 2 complex conjugate double roots & $D=0$ and $L1=0$ and\\ & $3p^{2}-8q<0$\\
8. 1 triple real root and 1 single real root & $D=0$ and $L1=0$ and\\ & $3p^{2}-8q>0$ and $p^{3}-4pq+8r\neq0$\\
8a. triple root $<$ single root & $p^{3}-4pq+8r<0$\\
8b. single root $<$ triple root & $p^{3}-4pq+8r>0$\\
9. quadruple root & $D=0$ and $L1=0$ and $3p^{2}-8q=0$\\
\end{tabular}\bigskip

\quad The above table gives two results.  If the lines labeled by pure numbers are chosen, then one gets the real and complex root multiplicities.  If all the lines are chosen, then one gets the real and complex root multiplicities together with the order of the real roots with respect to their multiplicities.  An automatic corollary of the above is the following: (Note that in 7.(Two complex conjugate double roots) of the real root configurations, the condition $p^{3}-4pq+8r=0$ holds, but was not needed there.  However, it is needed for the complex root configurations.)\newline

\textbf{Complex Root Configurations.}
\begin{tabular}{l r}
1. 4 distinct roots & $D\neq0$\\
2. 1 double root and 2 single roots & $D=0$ and $L1\neq0$\\
3. 2 double roots & $D=0$ and $L1=0$ and\\ & $3p^{2}-8q\neq0$ and $p^{3}-4pq+8r=0$\\
4. 1 triple root and 1 single root & $D=0$ and $L1=0$ and\\ & $3p^{2}-8q\neq0$ and $p^{3}-4pq+8r\neq0$\\
5. quadruple root & $D=0$ and $L1=0$ and $3p^{2}-8q=0$\\
\end{tabular}
\newpage

\section{Bibliography}
\noindent [1] Basu, Saugata, Pollack, Richard, and Roy, Marie-Francoise. \emph{Algorithms in Real Algebraic Geometry}, first edition. Springer-Verlag, 2003.\\
\noindent [2] Benedetti, Riccardo and Risler, Jean-Jacques. \emph{Real Algebraic and Semi-algebraic Sets}. Hermann, Editeurs des Sciences et des Arts, 1990.\\
\noindent [3] Jacobson, Nathan. \emph{Basic Algebra I}. 2nd ed. W.H. Freeman and Company. 1985.\\
\noindent [4] Liang, S. and Zhang. J. ``A complete discrimination system for polynomials with complex coefficients and its automatic generation.'' \emph{Science in China E} Vol. 42, No. 2, April 1999. p. 113-128.\\
\noindent [5] Sottile, Frank. \emph{Real Solutions to Equations from Geometry}. American Mathematical Society. 2011.\\
\noindent [6] Weinberg, David and Martin, Clyde. ``A note on resultants.'' \emph{Applied Mathematics and Computation} Vol. 24, 1987. p. 303-309.\\
\noindent [7] Yang, L., Hou, X.R., and Zeng, Z.B. ``A complete discrimination system for polynomials.'' \emph{Science in China E} Vol. 39, No. 6, Dec. 1996. p. 628-646.\\
\noindent [8] Yang, Lu. ``Recent advances on determining the number of real roots of paramateric polynomials.'' \emph{Journal of Symbolic Computation} Vol. 28, 1999. p. 225-242.\\

\end{document}